\documentclass[a4paper,11pt]{article}
\usepackage[latin1]{inputenc}
\usepackage[T1]{fontenc}
\usepackage[english]{babel}
\usepackage{amssymb}
\usepackage{amsfonts, amsmath}

\title{\textbf{The H\'{e}non-Heiles system as part of an integrable system
in five unknowns with three constants of motions}}
\author{\textbf{A. Lesfari}
\\\emph{Department of Mathematics}
\\\emph{Faculty of Sciences}
\\\emph{University of Choua\"{i}b Doukkali}
\\\emph{B.P. 20, El Jadida, Morocco}.
\\\emph{Lesfariahmed@yahoo.fr}}
\date{}
\begin{document}
\maketitle

\emph{Abstract}. In this paper we construct a new completely
integrable system. This system is an instance of a master system
of differential equations in $5$ unknowns having $3$ quartics
constants of motion.We find via the Painlevé analysis the
principal balances of the hamiltonian field defined by the
hamiltonian. Consequently, the system in question is algebraically
integrable. A careful analysis of this system reveals an intimate
rational relationship with a special case of the well known
H\'{e}non-Heiles system. The latter admits asymptotic solutions
with fractional powers in $t$ and depending on $3$ free
parameters. As a consequence, this system is algebraically completely integrable in the generalized sense. 

Among recent researches, there have been many works dealing with
the integrability of hamiltonian systems. In hamiltonian
mechanics, for integrating a dynamical system with $n$ degrees of
freedom, it is sufficient in most cases to know only the first $n$
integrals. This situation is known as Liouville complete
integrability of a hamiltonian system. A dynamical system is
algebraic completely integrable (in the sense of Adler-van
Moerbeke [1,2]) if it can be linearized on a complex algebraic
torus $\mathbb{C}^{n}/lattice$ (=abelian variety). The invariants
(often called first integrals or constants) of the motion are
polynomials and the phase space coordinates (or some algebraic
functions of these) restricted to a complex invariant variety
defined by putting these invariants equals to generic constants,
are meromorphic functions on an abelian variety. Moreover, in the
coordinates of this abelian variety, the flows (run with complex
time) generated by the constants of the motion are straight lines.
Also some interesting integrable systems appear as coverings of
algebraic completely integrable systems [4,9,12,14]. The invariant
varieties are coverings of abelian varieties and these systems are
called algebraic completely integrable in the generalized sense.
Beside the fact that many integrable hamiltonian systems have been
on the subject of powerful and beautiful theories of mathematics,
another motivation for its study is : the concepts of
integrability have been applied to an increasing number of
physical systems, biological phenomena, population dynamics,
chemical rate equations, to mention only a few. However, we have
no general method for testing the integrability of a given
hamiltonian system. It seems still hopless to describe, or even to
recognize with any facility, those hamiltonian systems which are
integrable, though they are quite exceptional. In the present
paper, we present an integrable system in five unknowns having two
cubic and one quartic invariants (constants of motion). Also, we
show that this system include in particular the well known
H\'{e}non-Heiles system.

Consider the system
\begin{eqnarray}
\stackrel{.}{y}_{1} &=&x_{1},  \nonumber \\
\stackrel{.}{y}_{2} &=&x_{2},  \\
\stackrel{.}{x}_{1} &=&-Ay_{1}-2y_{1}y_{2},  \nonumber \\
\stackrel{.}{x}_{2} &=&-By_{2}-y_{1}^{2}-\varepsilon y_{2}^{2},
\nonumber
\end{eqnarray}
corresponding to a generalized H\'{e}non-Heiles hamiltonian
$$
H=\frac{1}{2}(x_{1}^{2}+x_{2}^{2})+\frac{1}{2}%
(Ay_{1}^{2}+By_{2}^{2})+y_{1}^{2}y_{2}+\frac{\varepsilon
}{3}y_{2}^{3},
$$
where $A,$ $B,$ $\varepsilon $ are constant parameters and $
y_{1}$, $y_{2}$, $x_{1}$, $x_{2}$ are canonical coordinates and
momenta, respectively. First studied as a mathematical model to
describe the chaotic motion of a test star in an axisymmetric
galactic mean gravitational field [7], this system is widely
explored in other branches of physics. It well-known from
applications in stellar dynamics, statistical mechanics and
quantum mechanics$.$ It provides a model for the oscillations of
atoms in a three-atomic molecule [5].

Usually, the H\'{e}non-Heiles system is not integrable and
represents a classical example of chaotic behaviour. Nevertheless
at some special values of the parameters it is integrable; to be
precise, there are known three integrables cases :

$(i)$ $\varepsilon=6,$ $A$ and $B$ arbitrary. The second integral
of motion is
$$
H_{2}=y_{1}^{4}+4y_{1}^{2}y_{2}^{2}-4x_{1}^{2}y_{2}+4x_{1}x_{2}y_{1}+4Ay_{1}^{2}y_{2}+(4A-B)
x_{1}^{2}+A(4A-B)y_{1}^{2}.
$$

$(ii)$ $\varepsilon =1,$ $A=B$. The second integral of motion is
$$
H_{2}=x_{1}x_{2}+\frac{1}{3}y_{1}^{3}+y_{1}y_{2}^{2}+Ay_{1}y_{2}.
$$

$(iii)$ $\varepsilon =16,$ $B=16A.$ The second integral of motion
is
\begin{equation}\label{eqn:euler}
H_{2}=3x_{1}^{4}+6Ax_{1}^{2}y_{1}^{2}
+12x_{1}^{2}y_{1}^{2}y_{2}-4x_{1}x_{2}y_{1}^{3}-4Ay_{1}^{4}y_{2}-4y_{1}^{4}y_{2}^{2}+3A^{2}y_{1}^{4}-
\frac{2}{3}y_{1}^{6}.
\end{equation}

In the two cases (i) and (ii), the system (1) has been integrated
by making use of genus one and genus two theta functions. For the
case (i), it was shown [3] that this case separates in translated
parabolic coordinates. Solving the problem in case (ii) is not
difficult (this case trivially separates in cartesian
coordinates). In the case (iii), the system can also be integrated
[13] by making use of elliptic functions.

The general solutions of the equations of motion for hamiltonian
(1), for the case (i) and (ii), have the Painlev\'{e} propriety,
i.e., that they admit only poles in the complex time variable. For
the case (ii), we found [8,10,11] via Kowalewski-Painlev\'{e}
analysis the principal balances of the hamiltonian vector field
defined by the hamiltonian and we have shown that the system is
algebraic complete integrable. The affine surface defined by the
two constants of motion completes into an abelian surface  by
adjoining a smooth genus three hyperelliptic curve $D$. This curve
is a double cover of an elliptic curve $E$, ramified at four
points : it defines a line bundle and a polarization $(1,2)$ on
this abelian surface. Also, the abelian surface obtained is the
dual of a Prym variety, namely $\mbox{Prym}(D/E)^{\vee}$ and the
flow of solutions of the equations of motion is linearized on
abelian surface.

The present paper deals with the case $(iii)$. The system $(1)$
can be written in the form
\begin{equation}\label{eqn:euler}
\dot{u}=J\frac{\partial H}{\partial u},\quad
u=(y_{1},y_{2},x_{1},x_{2})^{\top},
\end{equation}
where
\begin{equation}\label{eqn:euler}
H\equiv H_1=\frac{1}{2}(x_{1}^{2}+x_{2}^{2})+\frac{A}{2}%
(y_{1}^{2}+16y_{2}^{2})+y_{1}^{2}y_{2}+\frac{16 }{3}y_{2}^{3},
\end{equation}
$$
\frac{\partial H}{\partial z}=(\frac{\partial H}{\partial y_{1}},\frac{%
\partial H}{\partial y_{2}},\frac{\partial H}{\partial x_{1}},\frac{\partial
H}{\partial x_{2}})^{\intercal },\mathrm{\ }\text{ }J=\left(
\begin{array}{ll}
O & I \\
-I & O
\end{array}
\right).
$$
The functions $H_1$(4) and $H_2$(2) commute :
$$\left\{ H_{1},H_{2}\right\} =\sum_{k=1}^{2}\left(
\frac{\partial H_{1}}{\partial p_{k}}\frac{\partial
H_{2}}{\partial q_{k}}-\frac{\partial H_{1}}{\partial
q_{k}}\frac{\partial H_{2}}{\partial p_{k}}\right)=0. $$ The
second flow commuting with the first is regulated by the equations
$$
\dot{u}=J\frac{\partial H_{2}}{\partial u},\quad
u=(y_{1},y_{2},x_{1},x_{2})^{\intercal }
$$
and is written explicitly as
\begin{eqnarray}
\dot{y}_{1} &=&-24Ax_{1}-8x_{1}y_{2}+4x_{2}y_{1},\nonumber\\
\dot{y}_{2} &=&4x_{1}y_{1},\nonumber\\
\dot{x}_{1} &=&24A^2y_{1}-\allowbreak
4x_{1}x_{2}-8Ay_{1}y_{2}-8y_{1}y_{2}^{2}-4y_{1}^{3},\nonumber\\
\dot{x}_{2} &=&4x_{1}^{2}-4Ay_{1}^{2}-8y_{1}^{2}y_{2},\nonumber
\end{eqnarray}

The invariant (or level) variety
\begin{equation}\label{eqn:euler}
A=\bigcap_{k=1}^{2}\{z\in\mathbb{C}^4:H_k(z)=b_k\},
\end{equation}
is a smooth affine surface for generic
$(b_{1},b_{2})\in\Bbb{C}^{2}$. It was conjectured in [13] that the
complex invariant variety $A$(5) is an affine part of an abelian
variety (an algebraic complex torus). We will see that in fact the
invariant surface $A$(5) can be completed as a cyclic double cover
$\overline{A}$ of an abelian surface $\widetilde{B}$. The system
(3) is algebraic complete integrable in the generalized sense.
Moreover, $\overline{A}$ is smooth except at the point lying over
the singularity of type $A_1$ : $X^2+Y^2-Z^2=0$ and the resolution
$\widetilde{A}$ of $\overline{A}$ is a surface of general type.

When one examines all possible singularities, one finds that it
possible for the variable $y_1$ to contain square root terms of
the type $t^{1/2}$, which are strictly not allowed by the Painlevé
test. However, these terms are trivially removed by introducing
some new variables $z_1,\ldots,z_5$ (see below), which restores
the Painlevé property to the system. The system (3) admits Laurent
solutions in $t^{1/2}$, depending on three free parameters:
$\alpha$, $\beta$, $\gamma$ and they are explicitly given as
follows
\begin{eqnarray}
y_{1}&=&\frac{\alpha }{\sqrt{t}}+\beta t\sqrt{t}-\frac{\alpha
}{18}t^{2}\sqrt{t}+\frac{\alpha
A^{2}}{10}t^{3}\sqrt{t}-\frac{\alpha ^{2}\beta}{18}t^{4}\sqrt{t}+\cdots,\nonumber\\
y_{2}&=&-\frac{3}{8t^{2}}-\frac{A}{2}+\frac{\alpha^{2}}{12}t-\frac{2A^{2}}{5}t^{2}
+\frac{\alpha\beta}{3}t^{3}-\gamma t^{4}+\cdots,\\
x_{1}&=&-\frac{1}{2}\frac{\alpha}{t\sqrt{t}}+\frac{3}{2}\beta\sqrt{t}-\frac{5}{36}\alpha
t\sqrt{t}+\frac{7}{20}\alpha A^{2}t^{2}\sqrt{t}-\frac{1}{4}\alpha^{2}\beta t^{3}\sqrt{t}+\cdots,\nonumber\\
x_{2}&=&\frac{3}{4t^{3}}+\frac{1}{12}\alpha^{2}-\frac{4}{5}A^{2}t+\alpha\beta
t^{2}-4\gamma t^{3}+\cdots\nonumber
\end{eqnarray}
These formal series solutions are convergent as a consequence of
the majorant method. By substituting these series in the constants
of the motion $H_{1}=b_{1}$ and $H_{2}=b_{2}$, i.e.,
\begin{eqnarray}
H_{1}&=&\frac{1}{9}\alpha ^{2}-\frac{21}{4}\gamma
+\frac{13}{288}\alpha ^{4}+\frac{4}{3}A^{3}=b_1,\nonumber\\
H_{2}&=&-144\alpha \beta ^{3}+\frac{294}{5}\alpha ^{3}\beta
A^{2}+\frac{8}{9}\alpha ^{6}-33\gamma \alpha ^{4}=b_2,\nonumber
\end{eqnarray}
one eliminates the parameter $\gamma$ linearly, leading to an
equation connecting the two remaining parameters $\alpha$ and
$\beta$ :
$$
144\alpha \beta ^{3}-\frac{294A^{2}}{5}\alpha ^{3}\beta
+\allowbreak \frac{143}{504}\alpha ^{8}-\frac{4}{21}\alpha
^{6}+\frac{44}{21}\left( 4A^{3}-3b_{1}\right) \alpha ^{4}+b_{2}=0.
$$
which is nothing but the equation of an algebraic curve
$\mathcal{D}$ of genus four along which the $u(t)
\equiv(y_{1}(t),y_{2}(t),x_{1}(t),x_{2}(t))$ blow up. To be more
precise $\mathcal{D}$ is the closure of the continuous components
of
$$
\left\{\text{Laurent series solutions }u(t)\text{ such that }
H_{k}(u(t))=b_{k},\text{ }1\leq k\leq 2\right\},
$$
i.e.,
$$
\mathcal{D}==t^{0}-\text{coefficient of }\left\{u\in
\Bbb{C}^{4}:\text{ }H_{1}(u(t))=b_{1}\right\}\cap\left\{u\in
\Bbb{C}^{4}:H_{2}(u(t))=b_{2}\right\} .
$$
The Laurent solutions restricted to the surface $A$(5) are
parameterized by this curve $\mathcal{D}$.

We show that the system (3) is part of a new system of
differential equations in five unknowns having two cubic and one
quartic invariants (constants of motion). By inspection of the
expansions (6), we look for polynomials in $(y_1,y_2,x_1,x_2)$
with a simple pole. Let
\begin{equation}\label{eqn:euler}
\varphi : A\longrightarrow \mathbb{C}^5,
(y_1,y_2,x_1,x_2)\longmapsto (z_1,z_2,z_3,z_4,z_5), \end{equation}
be a morphism on the affine variety $A(5)$ where $z_1,\ldots,z_5$
are defined as
$$z_1=y_{1}^{2},\quad z_2=y_{2},\quad z_3=x_{2},\quad
z_4=y_{1}x_{1},\quad z_5=3x_{1}^{2}+2y_{1}^{2}y_{2}.$$ The
morphism (7) maps the vector field (3) into a Liouville integrable
system (8) in five unknowns having two cubic and one quartic
invariants. Indeed, the change of variables (7) maps the vector
field (4) into the system of differential equation
\begin{eqnarray}
\dot{z}_1&=&2z_4,\nonumber\\
\dot{z}_2&=&z_3,\nonumber\\
\dot{z}_3&=&-z_{1}-16Az_{2}-16z_{2}^{2} ,\\
\dot{z}_4&=&-Az_{1}+\frac{1}{3}z_{5}-\frac{8}{3}z_{1}z_{2},\nonumber\\
\dot{z}_5&=&-6Az_{4}+2z_{1}z_{3}-8z_{2}z_{4},\nonumber
\end{eqnarray}
with constants of motion
\begin{eqnarray}
F_1&=&\frac{1}{2}Az_{1}+\frac{1}{6}z_{5}+8Az_{2}^{2}+\frac{1}{2}z_{3}^{2}
+\frac{2}{3}z_{1}z_{2}+\frac{16}{3}z_{2}^{3},\nonumber\\
F_2&=&9A^{2}z_{1}^{2}+z_{5}^{2}+6Az_{1}z_{5}-2z_{1}^{3}-24Az_{1}^{2}z_{2}
-\allowbreak 12z_{1}z_{3}z_{4}+24z_{2}z_{4}^{2}-16z_{1}^{2}z_{2}^{2},\nonumber\\
F_3&=&z_{1}z_{5}-3z_{4}^{2}-2z_{1}^{2}z_{2}.\nonumber
\end{eqnarray}
This new system is completely integrable and the hamiltonian
structure is defined by the Poisson bracket
$$\{F,H\}=\left\langle \frac{\partial F}{\partial z},
J\frac{\partial H}{\partial z}\right\rangle
=\sum_{k,l=1}^{5}J_{kl}\frac{\partial F}{\partial
z_{k}}\frac{\partial H}{\partial z_{l}},$$ where
$$\frac{\partial H}{\partial
z}=\left(\frac{\partial H}{\partial z_{1}},\frac{\partial
H}{\partial z_{2}},\frac{\partial H}{\partial
z_{3}},\frac{\partial H}{\partial z_{4}},\frac{\partial
H}{\partial z_{5}}\right)^\top,$$ and
$$J=\left(\begin{array}{ccccc}
0&0&0&2z_1&12z_4\\
0&0&1&0&0\\
0&-1&0&0&-2z_1\\
-2z_1&0&0&0&-8z_{1}z_{2}+2z_{5}\\
-12z_4&0&2z_1&8z_{1}z_{2}-2z_{5}&0
\end{array}\right),$$
is a skew-symmetric matrix for which the corresponding Poisson
bracket satisfies the Jacobi identities. The system (8) can be
written as
$$\dot{z}=J\frac{\partial
H}{\partial z},\quad z=(z_{1},z_{2},z_{3},z_{4},z_{5})^\top,$$
where $H=F_{1}.$ The second flow commuting with the first is
regulated by the equations
$$\dot{z}=J\frac{\partial F_{2}}{\partial z},\quad
z=(z_{1},z_{2},z_{3},z_{4},z_{5})^\top,$$ and is written
explicitly as
\begin{eqnarray}
\dot{z}_1&=&24z_{4}z_{5}-24z_{1}^{2}z_{3}+96z_{1}z_{2}z_{4}+72Az_{1}z_{4},\nonumber\\
\dot{z}_2&=&-12z_{1}z_{4},\nonumber\\
\dot{z}_3&=&12Az_{1}^{2}-12z_{1}z_{5}+48z_{1}^{2}z_{2},\nonumber\\
\dot{z}_4&=&4z_{5}^{2}-36A^{2}z_{1}^{2}+12z_{1}^{3}-16z_{1}z_{2}z_{5}+48Az_{1}^{2}z_{2}
+24z_{1}z_{3}z_{4}+64z_{1}^{2}z_{2}^{2},\nonumber\\
\dot{z}_5&=&-96z_{2}z_{4}z_{5}+768z_{1}z_{2}^{2}z_{4}-72Az_{4}z_{5}+576Az_{1}z_{2}z_{4}
+144z_{3}z_{4}^{2}\nonumber\\
&&\quad-216A^{2}z_{1}z_{4}+48z_{1}^{2}z_{4}-96z_{1}^{2}z_{2}z_{3}+24z_{1}z_{3}z_{5}.\nonumber
\end{eqnarray}
These vector fields are in involution, i.e.,
$$\{F_1,F_2\}=\left\langle \frac{\partial F_{1}}{\partial z},J\frac{\partial F_{2}}{\partial
z}\right\rangle=0,$$ and the remaining one is casimir, i.e.,
$$J\frac{\partial F_{3}}{\partial z}=0.$$
Consequently, the system (8) is integrable in the sense of
Liouville.

The invariant variety
\begin{equation}\label{eqn:euler}
B=\bigcap_{k=1}^{3}\{z\in\mathbb{C}^5:F_k(z)=c_k\},
\end{equation}
is a smooth affine surface for generic values of $c_{1},c_{2}$ and
$c_{3}$. The system (8) possesses Laurent series solutions which
depend on four free parameters. These meromorphic solutions
restricted to the surface B(9) can be read off from (6) and the
change of variable (7). Following the methods in [1,2,11], one
find the compactification of $B$ into an abelian surface
$\widetilde{B}$, the system of differential equations (8) is
algebraic complete integrable and the corresponding flows evolve
on $\widetilde{B}$. Applying the method explained in Piovan [20]
(see also [9]), we show that the invariant surface $A$(5) can be
completed as a cyclic double cover $\overline{A}$ of an abelian
surface $\widetilde{B}$. The system (3) is algebraic complete
integrable in the generalized sense. Moreover, $\overline{A}$ is
smooth except at the point lying over the singularity of type
$A_1$ : $X^2+Y^2-Z^2=0$ and the resolution $\widetilde{A}$ of
$\overline{A}$ is a surface of general type. We have shown that
the morphism $\varphi$ (9) maps the vector field (4) into an
algebraic completely integrable system (8) in five unknowns and
the affine variety $A$ (5) onto the affine part $B$ (9) of an
abelian variety $\widetilde{B}$. This explains (among other) why
the asymptotic solutions to the differential equations (3) contain
fractional powers.

In this paper we have shown that the Hénon-Heiles (third case)
equations are part of a system of differential equations in five
unknowns having three constants of motions. We have seen that this
new system is Liouville integrable but we do not know physical
interpretation for this system and it may be possible to have
links with other systems (in addition to Hénon-Heiles) with
applications in physics among other. The integrable dynamical
system presented here is an interesting example, particular to
experts of abelian varieties who may want to see an explicit
example of a correspondence for varieties defined by different
curves.

\end{document}